\documentclass[12pt]{amsart}

\def\thec{\textup{ceiling}(\log_2n)}

\def\answerRS{1+3\cdot\thec}
\def\ansRS{1+3c}
\def\answerR{2+3\cdot\thec}
\def\ansR{2+3c}

\def\answerIS{8+6\cdot\thec}
\def\ansIS{8+6c}
\def\answerI{9+6\cdot\thec}
\def\ansI{9+6c}

\usepackage{float}
\floatstyle{boxed} 
\restylefloat{figure}

\usepackage{amssymb}
\usepackage{amsmath}
\usepackage{color}

\begin{document}
\numberwithin{equation}{section} 

\def\1#1{\overline{#1}}
\def\2#1{\widetilde{#1}}
\def\3#1{\widehat{#1}}
\def\4#1{\mathbb{#1}}
\def\5#1{\mathfrak{#1}}
\def\6#1{{\mathcal{#1}}}
\def\mymf#1{\5#1}

\def\C{{\4C}}
\def\R{{\4R}}
\def\N{{\4N}}
\def\Z{{\4Z}}
\def\Q{{\4Q}}
\def\T{{\Theta}}

\def\sideremark#1{\ifvmode\leavevmode\fi\vadjust{
\vbox to0pt{\hbox to 0pt{\hskip\hsize\hskip1em
\vbox{\hsize1cm\tiny\raggedright\pretolerance10000
\noindent #1\hfill}\hss}\vbox to8pt{\vfil}\vss}}}
\def\sr#1{\sideremark{#1}}
\def\sr#1{}

\def\highlight#1{\color{red} #1\color{black}}
\def\hl#1{\highlight{#1}}
\def\hl#1{#1}

\title[]{Factoring Formal Maps into Reversible or Involutive Factors}
\author[A. O'Farrell \& D. Zaitsev]{Anthony G. O'Farrell
\& Dmitri Zaitsev}
\address{%
A. O'Farrell: Department of Mathematics and Statistics,
 NUI, Maynooth, Co. Kildare, Ireland}
 \email{anthonyg.ofarrell@gmail.com}
 \address{
D. Zaitsev: School of Mathematics, 
Trinity College Dublin, Dublin 2, Ireland
}
\email{zaitsev@maths.tcd.ie}

\thanks{
Anthony G. O'Farrell: Mathematics and Statistics, NUI, Maynooth, Co. Kildare, Ireland.\\ 
e-mail: anthonyg.ofarrell@gmail.com\\
Dmitri Zaitsev: Mathematics, Trinity College, Dublin 2, Ireland.\\ 
email: zaitsev@maths.tcd.ie\\
Supported in part by the Science Foundation Ireland grant 10/RFP/MTH2878.
}

\def\Label#1{\label{#1}}


\def\cn{{\C^n}}
\def\cnn{{\C^{n'}}}
\def\ocn{\2{\C^n}}
\def\ocnn{\2{\C^{n'}}}


\def\dist{{\rm dist}}
\def\diag{{\rm diag}}
\def\const{{\rm const}}
\def\rk{{\rm rank\,}}
\def\id{{\sf id}}
\def\aut{{\sf aut}}
\def\Aut{{\sf Aut}}
\def\CR{{\rm CR}}
\def\GL{{\sf GL}}
\def\SL{{\sf SL}}
\def\U{{\sf U}}
\def\Re{{\sf Re}\,}
\def\Im{{\sf Im}\,}
\def\im{{\rm im}\,}
\def\span{\text{\rm span}}
\def\tr{{\sf tr}\,}
\def\HOT{{\textup{\small HOT}}}

\def\codim{{\rm codim}}
\def\crd{\dim_{{\rm CR}}}
\def\crc{{\rm codim_{CR}}}

\def\phi{\varphi}
\def\eps{\varepsilon}
\def\d{\partial}
\def\a{\alpha}
\def\b{\beta}
\def\g{\gamma}
\def\G{\Gamma}
\def\D{\Delta}
\def\Om{\Omega}
\def\k{\kappa}
\def\l{\lambda}
\def\L{\Lambda}
\def\z{{\bar z}}
\def\w{{\bar w}}
\def\Z{{\mathbb Z}}
\def\t{\tau}
\def\th{\theta}
\def\Oh{\textup{O}}

\emergencystretch15pt
\frenchspacing

\newtheorem{Thm}{Theorem}[section]
\newtheorem{Cor}[Thm]{Corollary}
\newtheorem{Pro}[Thm]{Proposition}
\newtheorem{Lem}[Thm]{Lemma}

\theoremstyle{definition}\newtheorem{Def}[Thm]{Definition}

\theoremstyle{remark}
\newtheorem{Rem}[Thm]{Remark}
\newtheorem{Exa}[Thm]{Example}
\newtheorem{Exs}[Thm]{Examples}

\def\bl{\begin{Lem}}
\def\el{\end{Lem}}
\def\bp{\begin{Pro}}
\def\ep{\end{Pro}}
\def\bt{\begin{Thm}}
\def\et{\end{Thm}}
\def\bc{\begin{Cor}}
\def\ec{\end{Cor}}
\def\bd{\begin{Def}}
\def\ed{\end{Def}}
\def\br{\begin{Rem}}
\def\er{\end{Rem}}
\def\be{\begin{Exa}}
\def\ee{\end{Exa}}
\def\bpf{\begin{proof}}
\def\epf{\end{proof}}
\def\ben{\begin{enumerate}}
\def\een{\end{enumerate}}
\def\beq{\begin{equation}}
\def\eeq{\end{equation}}

\def\dsty{\displaystyle }

\begin{abstract}
An element $g$ of a group is called {\em reversible}
if it is conjugate in the group to its inverse.  
An element is an {\em involution}
if it is equal to its inverse.
This paper is about factoring elements
as products of reversibles in the group
$\5G_n$ of formal maps
of $(\C^n,0)$, i.e. formally-invertible
$n$-tuples of formal power series in $n$
variables, with complex coefficients.
The case $n=1$ was already understood \cite{O}.

Each product $F$ of reversibles has linear part $L(F)$ of determinant
$\pm1$. The main results are that
\sr{absorbed $n=2$ case}
for $n\ge2$ each map $F$ with $\det(L(F))=\pm1$ 
is the product of $\answerR$
reversibles, and may
also be factored
as the product of $\answerI$
involutions
(where the ceiling of $x$ is the smallest
integer $\ge x$).
\end{abstract}

\maketitle

\section{Introduction}
\subsection{}
It is an interesting fact that in many very large
groups each element may be factored as the product of 
a small number of involutions.
For instance, each permutation is the product of two involutions.
Less trivially, Fine and Schweigert \cite{FS} showed
that each homeomorphism of $\mathbb{R}$ onto itself 
is the composition at most four involutions, each one conjugate to
the map $x\mapsto -x$. 

A natural generalization of an involution is a reversible.
An element $g$ of a group is called {\em reversible}
if it is conjugate to its inverse, i.e. the conjugate $g^h:=h^{-1}gh$ equals $g^{-1}$ for some $h$ from the group.  We say that $h$ {\em reverses} $g$ or $h$ {\em is a  reverser of} $g$, in this case.
Furthermore, if the reverser $h$ can be chosen to be an involution
(i.e. an element of order at most $2$), then $g$ is called {\em strongly reversible}.
(Note that some writers use the terminology ``weakly reversible'' and ``reversible''
instead of respectively ``reversible'' and ``strongly reversible'' used here. In finite group theory,
the terms used are \lq\lq real" and \lq\lq strongly real".)
A strongly reversible element is the product of two involutions.
See \cite{Go76, Go79, Go81}.
If $g$ is reversed by an element of finite even order $2k$,
then $g$ is the product of two elements of order $2k$.
\sr{added text}
Indeed, it is easy to check that if $g$ is reversed by
some element $h$, then it factors as $hf$, where 
$h^2=f^{-2}$, so if $h$ has order $2k$, then
so does $f$.

Reversible maps have their origin in problems of classical dynamics,
such as the harmonic oscilator, the $n$-body problem or 
billiards,
and Birkhoff \cite{B} was one of the first to realize their significance.
He observed that a Hamiltonian system with Hamiltonian quadratic in 
the momentum (such as the $n$-body problem), and other
interesting dynamical systems 
admit what are called ``time reversal symmetries'',
i.e. transformations of the phase space that
conjugate the dynamical system to its inverse. 

In CR geometry reversible maps played important role in the celebrated 
work of Moser and Webster \cite{MW},
arising as products of two involutions
naturally associated to a CR singularity.
Such a reversible map is called there ``a discrete version of the Levi form"
and plays a fundamental role
in the proof of the convergence of the normal form
for a CR singularity.
More recently, this map has been used
by Ahern and Gong \cite{AG}
for so-called parabolic CR singularities.

The basic concept of
reversible element makes sense in any group,
and reversibility has been the focus of interest
in many other application areas that 
involve some underlying group.
For instance, reversible elements
appear (sometimes under aliases) in
connection with geometrical symmetries,
special geodesics on
Riemann surfaces, binary integral quadratic forms, 
quadratic correspondences, superposition of functions,
approximation problems, toral automorphisms and foliations
\cite{Dj1967, Dj1986, FZ, GS, Ni1987, Sar}.
For further references to some contexts in which reversible elements
have played a part, and a short survey of factorization
results involving reversibles, see \cite{OCRM}.

\hl{From the point of view of group theory, 
the subgroup $R^\infty(G)$ generated by the reversible 
elements of a group $G$ is normal, and 
its isomorphism class is an isomorphism
invariant of $G$. It has associated
numerical invariants, which are very basic
invariants of $G$, and their determination
is a natural first step in the classification of
$G$.  One of these invariants is the 
supremum over all $g\in R^\infty(G)$ of the
\sr{rewritten}
least number $k$ of reversible factors $r_j$ 
needed to represent $g$
as a product $r_1\cdots r_k$. In the language 
of Klopsch and Lev \cite{KL}, this is
the \lq\lq diameter'' of $R^\infty$ with respect to
the set $R(G)$ of reversibles. 
}

\sr{added citations}
\hl{
The issue of factorization into reversibles (and
involutions), and the number of factors needed
 has attracted attention in several 
group contexts --- see for instance 
\cite{BrFa, El1993, GHR, Is, KnNi, St}.
}

In this paper we consider the group $\5G_n$
of formally-invertible maps in $n\ge2$ 
complex variables, and we 
discuss the factorization of a given map
as a product of reversibles, and as a product
of involutions. \hl{We get an explicit upper bound
in terms of $n$ for the above diameter,
and also (when $n\ge2$) for the 
(finite!) diameter of $R^\infty(\5G_n)$
with respect to the set of involutions. 
}

In previous work
the first author dealt with this problem 
for $n=1$, and obtained the following results:

\bt\cite{O}\Label{T:O}
Let $F\in\5G_1$. 
Then the following are equivalent:\\
(1) $F$ is a product of reversibles.
\\
(2) $F(z)=\pm z+$ $\Oh(z^2)$,
\hl{i.e. $F(z) = \pm z + $terms in $z^2$ and higher
powers of $z$.}
\\
(3) $F$ is the product of two reversibles.
\qed
\et

\bt\cite{O}\Label{T:O-I}
Let $F\in\5G_1$. Then  the following are equivalent:
\\
(1) $F$ is a product of involutions.
\\
(2) For some $a\in\C$,
$F(z)=\pm z+ az^2 \pm a^2z^3+\Oh(z^4)$.
\\
(3)
$F$ is the product of four involutions.
\qed
\et

Thus not every reversible series in one variable
is the product of a finite number of involutions.  It depends upon the
conjugacy class of the series, modulo $z^4$.
We shall see that the
 situation changes in higher dimensions.

In dimension $2$, the authors previously showed
the following:
\sr{reduced to old result only}
\bt\cite{OZ}\Label{T:prod-4-rev}
If $F\in\5G_2$ has linear part of determinant $1$, then it may be factorized
\\
as the product
of $4$ reversible elements.
\et


\subsection{Results}
In this paper, we will show:
\sr{folded in $n=2$}
\bt\Label{T:main}
Let $n\ge 2$ and $F\in\5G_n$
have linear part of determinant $1$.
Let $c=\thec$.
Then\\ 
(1) $F$ is the product of $\ansRS$ reversibles.
\\
(2) $F$ is the product of $\ansIS$ involutions.
\qed
\et 

We also have:

\bc\Label{C:main}
Let $n\ge 2$ and $F\in\5G_n$. Let
$c=\thec$.
Then the following are equivalent:\\
(1) $F$ is a product of reversibles.
\\
(2) The linear part of $F$ has determinant $\pm1$.
\\
(3) $F$ is the product of $\ansR$ reversibles.
\\
(4) $F$ is the product of $\ansI$ involutions.\\
(5) 
$F$ is the product of $3+6c$ involutions
and two reversible maps of order dividing $4$.\sr{changed degree to order}
\qed
\ec
 
Thus, for instance in dimension 2, {\em every}
product of reversibles is also the product of
at most 15 involutions.

\subsection{Outline}
In Section \ref{S:Prelim} we define terminology
and notation, and develop some tools
that will be used in the proofs of these results.
We identify some interesting subgroups of
$\5G_n$, and construct homomorphisms 
connecting them.  In particular, we identify
a subgroup $\5C_n$, the centraliser in $\5G_n$ of
a matrix subgroup $D_n\le\GL(n,\C)\le\5G_n$, and we represent
$\5C_n$ as the semidirect product
of an abelian subgroup all of whose elements
are reversible in $\5G_n$ and a subgroup
(called $\5K_k$ or $\hat{\5K_k}$,
depending on whether $n$ is 
even or odd) of $\5C_k$, where $k$ is roughly
half of $n$.  This structural information
\sr{reference to exact sequences}
is summarized in the exact sequences
shown in Figure 
\ref{fig:exact-sequence2} below.
This allows us to
carry out an induction, reducing the 
reversible factorization of elements 
of $\5C_n$ to the reversible factorization of
$k$-dimensional maps, at the cost of
one extra factor.  Also,
the subgroup $\5C_n$ has a representative
of each so-called {\em generic}
conjugacy class in $\5G_n$,
and at the cost of an extra couple of factors,
we can reduce the factorization
of a general element of $\5G_n$ to
the factorization of a generic element.

These subgroups and homomorphisms elaborate upon
tools that were employed in our previous paper
\cite{OZ}, in which we characterized the
generic reversibles in dimension $2$.

In considering involutive factors, we have to deal with
the fact that not all one-dimensional maps $\chi\in\5G_1$
with multiplier $1$ can be factored into involutions,
so we have to find a way to
factor the lift $H(\chi)\in\5G_2$
into involutions.  Once we manage to do this,
we can then start the induction at $n=2$
and continue as before.  This depends on the fact that
the extra two or three reversible factors needed
at each induction step are all {\em strongly reversible},
i.e. products of two involutions (see below).  

In Section \ref{S:D-2}
we prove the two-dimensional results,
and in Section \ref{S:D-3}
we prove the rest.
 
\subsection{Open Questions}\sr{added}
When we get into the detailed proofs,
it will appear that for certain
dimensions $n$ we can derive
estimates for the number of reversible factors needed
that are considerably
smaller than the estimate in Theorem \ref{T:main}.
For instance,
we can do much better with $n=96$
than $n=97$.
See Section \ref{S:D-3} and Table \ref{Table} for details.

But we do not know sharp values for the number
of reversible or involutive factors needed
in any case
of dimension greater than $1$. It may even be the case
that a universal number of factors suffices
in all dimensions. 
Also,
it remains open, even for one-variable maps, whether
results such as these hold for convergent power series.
These are
interesting problems.  

\hl{One might wonder whether the coefficient field
$\C$ may be replaced by another in these results.
In our arguments, the properties of $\C$
that we use are the fact that it has characteristic
zero and  is algebraically-closed.
\sr{bit about rank removed}
We have not investigated
more general fields.  The paper \cite{O} gave
a complete account of reversibility and factorization
into reversibles in the one-dimensional 
formal map group for arbitrary coefficient fields
of characteristic zero.  As far as we know, there
is little known about reversibility 
when the characteristic of the coefficient field
is finite.  One should mention
that, thanks to Klopsch \cite[p.16]{C}, \cite{K}
the involutions (and indeed 
the elements of finite order)
have been identified  
for the so-called Nottingham groups (the one-dimensional
case in which 
the coefficient field is finite),
at least when the order 
of the field is odd.
}
 
\section{Notation and Preliminaries}\Label{S:Prelim}
\subsection{Power Series Structures}
For $n\in\N$,
let $\5F_n$ denote the ring of 
formal power series in $n$ 
(commuting) variables, with
complex coefficients,
and let $\5F_n^\times$ denote
the multiplicative group
of
its invertible elements,
i.e.\ those with nonzero constant term,
and let $\5M_n$ denote the complementary set
$\5F_n\sim\5F_n^\times$, the maximal ideal. 
Then an element of the set $\5S_n=(\5M_n)^n$
of $n$-tuples of elements of $\5M_n$
may be thought of as a formal map
of $\C^n\to\C^n$, taking $0$ to $0$.
Under formal composition,  $\5S_n$
is a semigroup, with identity
$\id_n(z)=z$.
Let $\5G_n\subset\5S_n $
be the group $\5G_n$
of formally-invertible elements.

We remark that $\5G_n$ is
isomorphic to the group of $\C$-algebra 
automorphisms of $\5F_n$. 
Indeed,
if $g\in\5G_n$, then
$f\mapsto f\circ g$ is
an automorphism of $\5F_n$.
Conversely, let $\Phi$ be any automorphism
of $\5F_n$, and take $g=(\Phi(z_1),\cdots,
\Phi(z_n))$.  Then 
$\Phi$
must map the unique maximal ideal
$\5M_n$ onto itself, and hence 
determines an automorphism of each
quotient $\5F_n/(\5M_n^k)$.
Since (the cosets of) $z_1$,$\ldots$,$z_n$
generate    $\5F_n/(\5M_n^k)$,
we have $\Phi(f)=f\circ g$ mod
$\5M_n^k$ for each $f$. Since this holds for each
$k\in\N$, we conclude that
$\Phi$ is just $f\mapsto f\circ g$.
\sr{removed Frechet, etc.}

\subsection{The map $L$}\label{l-notation}
A typical element $F\in\5S_n$ takes the form
$$F(z)=(F_{,1}(z),\ldots,F_{,n}(z))=(F_{,1}(z_1,\ldots,z_n),\ldots,
F_{,n}(z_1,\ldots,z_n))$$
where each $F_{,j}(z)$ is a power series in $n$ variables
having complex coefficients, and no constant term.
We shall refer to such series $F$ as maps, even though
they may be just \lq formal', i.e. 
the series may fail to converge at any $z\not=0$.

We usually write the formal composition of two maps
$F,G\in\5S_n$ as $FG$. We also write the product of two complex
numbers $a$ and $b$ as $ab$, but in cases where there might
be some ambiguity we use $a\cdot b$. For $n$-tuples
$a=(a_1,\ldots,a_n)$ and $b=(b_1,\ldots,b_n)$ (of various kinds)
we also use $a\cdot b$ for the \lq dot' product
$a_1\cdot b_1+\cdots+a_n\cdot b_n$,
and, a little more unusually, we will use
$a\times b$ for the {\em coordinatewise product}:
$$ (a_1,\ldots,a_n)\times
(b_1,\ldots,b_n):= (a_1b_1,\ldots,a_nb_n).$$ 

\smallskip
The series $F$ may be expressed
as a sum  
$$F = \sum_{k=1}^\infty L_k(F), $$
where $L_k(F)$ is homogeneous of degree $k$. 
We abbreviate $L_1(F)$ to $L(F)$. This term,
the {\em linear part of} $F$, belongs to
the algebra of $n\times n$ matrices.

An element $F$ of $\5S_n$
belongs to $\5G_n$ if and only its
linear part $L(F)$ belongs to 
the general linear group $\GL(n,\C)$.

We have the inclusion
$\GL(n,\C) \to \5G_n$,
and $L:\5G_n\to \GL(n,\C)$ is a group homomorphism.
We always identify $\GL(n,\C)$
with its image in $\5G_n$.

The elements of the kernel of $L$ 
are said to be {\em tangent to the identity}.

\subsection{Elements of Finite Order}
We note the following \cite[Lemma 2.1]{OZ}:
\begin{Lem}\Label{L:fol}
Let $n\in\N$ and let $\5H$ be a subgroup of $\5G_n$
such that\\
(1) $L(F)\in \5H$ whenever $F\in\5H$, and
\\
(2) $\5H\cap\ker L$ is closed under
convex combinations, i.e.
if $F_1,F_2\in\5H$, $L(F_1)=L(F_2)=\id$
and $0<\alpha<1$, then $\alpha F_1+(1-\alpha)F_2\in\5H$.
\\
Suppose $\T\in\5H$ has finite order. Then $\T$ is conjugated
by an element of $\5H\cap \ker L$ to its linear part $L(\T)$.
\end{Lem} 

This applies to $\5H=\5G_n$, $\5G_n\cap\ker L$,
 $\5G_n\cap\ker(\det\circ L)=L^{-1}(\SL(n,\C))$
\sr{removed unitary gp}
 (and, more generally to $L^{-1}(H)$
for any subgroup $H\le\GL(n,\C)$), to the corresponding
subgroups of biholomorphic germs (i.e. series that converge
on a neighbourhood of the origin) and to other subgroups introduced
below.  It applies to the intersection of any two groups
to which it applies.

In particular, in any $\5H$ to which the lemma applies,
each involution is conjugate to one of the linear involutions
in the group.  In $\GL(n,\C)$, a matrix is an involution
if and only if it is diagonalizable with eigenvalues $\pm1$.

Thus the involutions in $\5G_n$ are all conjugate
to their linear parts, which are 
involutions in $\GL(n,\C)$, and are classified up to
conjugacy by the dimension of the eigenspace
of the eigenvalue $1$. Thus there are just
$n$ conjugacy classes of proper involutions,
and condition (4) in Corollary \ref{C:main} says that
for $n\ge3$ one may represent any such $F$ as the product
of at most  $\answerI$ elements drawn
from this small collection of classes.

We remark that there are also just a finite number
of conjugacy classes in $\5G_n$ of maps of order dividing $4$.
The number is the number of ordered partitions of
$n$ as a sum of $4$ nonnegative integers, which 
equals $\binom{n+3}3$.

\subsection{Linear reversibles}\label{linear}
Reversibility is preserved by homomorphisms,
so a map $F\in\5G_{n}$ is reversible only if 
$L(F)$ is reversible in $\GL(n,\C)$. 
Classification of linear reversible maps is simple.
Suppose $F\in\GL(n,\C)$ is reversible.
Since the Jordan normal form of $F^{-1}$
consists of blocks of the same size as $F$
with inverse eigenvalues, 
the eigenvalues of $F$ that are not $\pm1$
must split into groups of pairs $\l,\l^{-1}$.
Furthermore, we must have the same number of Jordan blocks
of each size for $\l$ as for $\l^{-1}$.
Vice versa, if the eigenvalues of $F$ are either $\pm1$
or split into groups of pairs $\l,\l^{-1}$
with the same number of Jordan blocks of each size,
then both $F$ and $F^{-1}$ have the same Jordan normal form
and are therefore conjugate to each other.

\subsection{The Groups $D\le\GL(2,\C)$ and $D_n\in\GL(n,\C)$}
In particular, a linear map
is reversible in $\GL(2,\C)$
if and only if 
it is an involution or is conjugate to 
$\left(\begin{matrix}1&1\\0&1\end{matrix}\right)$,
$\left(\begin{matrix}-1&1\\0&-1\end{matrix}\right)$ 
or to
a matrix of the form
\beq\Label{E:glg}
 \tau_\mu=\left(\begin{matrix}\mu&0\\0&
\displaystyle\mu^{-1}\end{matrix}\right),
\eeq
for some $\mu\in\C^{\times}$.  
Thus each reversible $F\in\5G_2$ is conjugate \sr{added subscript 2}
in $\5G_2$ (by a linear conjugacy) to
a map having one of these types as its linear part.

The collection of maps $\tau_\mu$, defined in \eqref{E:glg}
forms an abelian subgroup of $\5G_2$, which we denoted by $D$
in \cite{OZ}.
\sr{changed tense. Added citation}
The element \eqref{E:glg} has infinite order
precisely when $\mu$ is not a root of unity,
and this is what we regarded as the generic situation
when $n=2$.

\medskip
We now extend this notation to higher dimensions.

\sr{merged subsections}
When $n=2m\ge2$ is even, we denote by $D_n$
the set of maps $T\in\5G_n$ of the form
$$ T(z) = \left(T_1(z_1,z_2),\ldots,T_m(z_{n-1},z_n)\right),$$
where each $T_j\in D$.

When $n=2m+1\ge3$ is odd, we denote by $D_n$
the set of maps $T\in\5G_n$ of the form
$$ T(z) =  \left(T_1(z_1,z_2),\ldots,T_m(z_{n-2},z_{n-1}),z_n\right),$$
where each $T_j\in D$, i.e. 
$T = T'\times\id_1$, where $T'\in D_{n-1}$ and $\id_1$ is the identity
map of $\C$.

In either case ($n=2m$ or $n=2m+1$), 
$D_n$ is a subgroup of $\5G_n$,
isomorphic to the $m$-fold cartesian product $D^m$.

An element $T\in D_n$ is called {\em generic} if
the associated $T_j=\tau(\mu_j)$, where there is no
\lq\lq resonance" relation
$$ \mu_1^{r_1}\cdots\mu_m^{r_m}=1 $$
with each $r_j\in\Z$, except the trivial
relation with all $r_j=0$.  If $T$ is generic, then
in particular no $\mu_j$ is a root of unity.
One could rephrase the condition as stating
that the $\mu_j$ generate a free abelian subgroup
of $\C^\times$ of rank $m$.

\medskip
We shall make use of the
classical Poincar\'e-Dulac
Theorem \cite[Section 4.8, Theorem 4.22]{PD},
and we state it here in our language,
for the reader's convenience:

\bt[Poincar\'e-Dulac]
Each map $F\in\5G_{n}$
is conjugate in $\5G_n$ 
to a map in the centralizer in $\5G_n$ 
of the linear part $L(F)$.
\qed\et

In case $L(F)$ is a generic member of $D_n$
we shall see shortly  (cf. Lemma \ref{L:Centraliser})
that the centralizer of $L(F)$ in $\5G_n$
coincides with the centralizer of the whole
subgroup $D_n$ in $\5G_n$.

\subsection{The Group $\5C_n=C_{D_n}(\5G_n)$}\label{cent-section}
In what follows, we shall usually have to distinguish
odd and even $n\ge2$.  When $z\in\C^n$ with 
$n=2m$ or $n=2m+1$, we define
$$ p(z) = (z_1z_2,\ldots,z_{2m-1}z_{2m}), $$
and we set
$$\pi(z):=
\left\{
\begin{array}{rcl}
p &,& n=2m,\\
(p,z_n) &,& n=2m+1,
\end{array}
\right.
$$
Both $p$ and $\pi$ depend (implicitly) on $n$.  

It is convenient, when dealing with 
$\5G_n$ for a given $n\ge2$,
to denote by $k$
the number 
$$
k=\left\{
\begin{array}{rcl}
m &,& n=2m,\\
m+1 &,& n=2m+1.
\end{array}
\right.
$$
Thus $m$ is the floor of $n/2$,
and $k$ is its ceiling.
We shall assume this relation between
$n$, $m$ and $k$, always.

The map $\pi$ sends $\C^n$ onto $\C^k$. A right
inverse for $\pi$ is the map $\epsilon=\epsilon_n \colon \C^k\to\C^n$,
given by
$$
\epsilon(t_1,\ldots,t_m)=(t_1,1,\ldots,t_m,1)$$
when $n=2m$,
and
$$
\epsilon(t_1,\ldots,t_m,t_k)=(t_1,1,\ldots,t_m,1,t_k)$$
when $n=2m+1$.

We note, for future reference, that
these maps preserve the coordinatewise product:
$$ \pi( a\times b ) = \pi(a)\times\pi(b),\ \forall a,b\in\C^n$$
and
$$ \epsilon( a\times b ) = \epsilon(a)\times\epsilon(b),\ \forall a,b\in\C^k.$$

\bl\Label{L:Centraliser}
Let $n\ge2$, and $F\in\5G_n$. 
Then the following are equivalent:\\
(1) $F$ commutes with each element of $D_n$.
\\
(2) For some generic $\L\in D_n$, $F$ commutes with $\L$. 
\\
(3)
\\
If $n=2m$ is even, then
$F$ takes the form 
$$F(z) = \left(
z_1\phi_1(p),
\ldots,
z_{n}\phi_{n}(p)
\right),$$
for some $\phi_j\in\5F_m^\times$.
\\
If $n=2m+1$ is odd, 
then
$F$ takes the form 
$$F(z) = \left(
z_1\phi_1(p,z_n),
\ldots,
z_{2m}\phi_{2m}(p,z_n),
z_{n}\phi_n(p,z_n)+p_1\zeta_1(p)+\cdots+p_m\zeta_m(p)
\right)
$$
for some $\phi_j\in\5F_{m+1}^\times$ and some $\zeta_k\in\5F_m$.
\el

\bpf
The only nontrivial implication is (2)$\implies$(3).

Suppose (2) holds, and fix $F\in\5G_n$ commuting with
some generic $\L\in D_n$.

\medskip
Case $1^\circ$: $n=2m$, even.

Then $\L=\diag(\l_1,1/\l_1,\ldots,\l_m,1/\l_m)$,
and there is no nontrivial resonance relation
$\prod \l_j^{\k_j}=1$.

We may write $F_{,1}$, the first component of $F$,
in the form
\sr{changed all $\phi$ to $\psi$} 
$$ z_1\psi_1(z_1,\ldots,z_n)
+
z_2\psi_2(z_2,\ldots,z_n)
+
\cdots+
z_n\psi_n(z_n),$$
where $\psi_j\in\5F_{n-j+1}$
(by gathering all monomial terms that involve $z_1$
into the term 
$z_1\psi_1(z_1,\ldots,z_n)$, all the terms that
involve $z_2$ but not $z_1$ into the next,
and so on).  Equating the first components
in the two sides of the equation $F\L=\L F$
gives
$$\begin{array}{rcl}
&&
\l_1z_1\psi_1(\L z)
+\dsty\frac{z_2}{\l_1}\psi_2\left(
\frac{z_2}{\l_1},\l_2z_3,\ldots\right)
+\cdots+
\frac{z_n}{\l_m}\psi_n\left(
\frac{z_n}{\l_m}\right)
\\
&=&
\l_1 z_1\psi_1(z_1,\ldots,z_n)
+
\l_1z_2\psi_2+
\cdots+
\l_1z_n\psi_n(z_n).
\end{array}
$$
Now, equating the coefficients
of each monomial on the two sides, and using
nonresonance, gives
\sr{expand??}
that 
$\psi_2=\cdots=\psi_n=0$
and
$\psi_1(\L z) = \psi_1(z)$, so that
$\psi_1(z)$ depends only on $p$.
Thus the first component of $F$ has the desired form.

A similar argument shows that each other component
takes the form in 
(3), so (3) holds.

\medskip
Case $2^\circ$: $n=2m+1$, odd. 
This time $\L=\diag(\l_1,1/\l_1,\ldots,\l_m,1/\l_m,1)$,
and again there is no nontrivial resonance relation
between the $\l_j$, $j=1$,$\ldots$,$m$.

Focussing, as before, on the first component
in the identity $F\L=\L F$, 
we have
$$\begin{array}{rcl}
&&
\l_1z_1\psi(\L z)
+\dsty\frac{z_2}{\l_1}\psi_2\left(
\frac{z_2}{\l_1},\l_2z_3,\ldots\right)
+\cdots+
\frac{z_{2m}}{\l_m}\psi_n\left(
\frac{z_{2m}}{\l_m},z_n\right)
+z_n\psi_n(z_n)
\\
&=&
\l_1 z_1\psi_1(z_1,\ldots,z_n)
+
\l_1z_2\psi_2+
\cdots+
\l_1z_{2m}\psi_{2m}(z_{2m},z_n)
+
\l_1z_n\psi_n(z_n).
\end{array}
$$
Identifying terms, as before,
we see that it proceeds just as in the even case
(with $z_n$ as an added parameter),
for $\psi_1$,$\ldots$,$\psi_{2m}$,
and find that $\psi_1$ depends only
on $p$ and $z_n$, and that $\psi_2=\cdots=\psi_{2m}=0$.
Finally, 
$\psi_n(z_n)=0$, since $\l_1\not=1$,
so that the first component of $F$ takes the
desired form.

A similar argument looks after all the components
except the last.  

Writing $z=(z',z_n)$, with $z'\in\C^{2m}$, we 
may write the $n$-th component of $F$ in the form
$$ F_{,n}(z) = 
z'\cdot G(z') + z_n\psi_n(z),$$
where $\psi_n\in\5F_n$ and $G(z')\in(\5F_{2m})^{2m}$ 
is a $2m$-vector of power series
in $2m$ variables, and $\cdot$ here denotes 
the dot product. (This $\psi_n$ is not
the one used in the argument about
the first component, the one that turned out to
be zero.) 
The last component of the identity $F\L=\L F$
then yields
$$
(\L'z')\cdot G(\L'z') + z_n\psi_n(\L z)
= z'\cdot G(z') + z_n\psi_n(z),  
$$
where $\L'$ denotes 
$\diag(\l_1,1/\l_1,\ldots,\l_m,1/\l_m)$.
This tells us that $\psi_n(z)$ depends
only on $p$ and $z_n$, and that
$z'\cdot G(z')$ depends only on $p$,
and hence takes the form $p\cdot\zeta(p)$,
for some $m$-tuple $\zeta\in(\5F_m)^m$.
Thus (3) holds. 
\epf

\br
We note that by condition (3) of the lemma,
each $F\in\5C_n$ has a diagonal linear part $L(F)$, because the terms
$p\cdot\zeta(p)$ that occur in the odd case 
are at least quadratic, so that in all cases
$L(F)=\diag(\phi_1(0),\ldots,\phi_n(0))$.
\er

\bd
We denote by $\5C_n$ the group of all maps
$F\in\5G_n$ that satisfy
any of the equivalent conditions of Lemma  \ref{L:Centraliser}.
\ed

\subsection{The Functions $M$, $\hat M$ and the Involution $J$} \sr{changed subsection title}
In terms of the coordinatewise product,
in the even case $n=2m$
we may represent the $F(z)$ in condition (3) 
of the lemma more compactly
as $z\times\phi(p)$.  We also denote this
map $F$ by $M(\phi)$.  Thus $M=M_n$ is a bijection
from $(\5F_m^\times)^{2m}$ onto
$\5C_{2m}$.  It is {\em not}, however, a homomorphism
from the abelian product group structure
of $(\5F_m^\times)^{2m}$.

\smallskip
For $n=2m+1$ odd, we note from the proof of the lemma
that
for $F\in\5C_n$, the last component $F(z)_{,n}$
takes the form $\psi(\pi(z))=\psi(p,z_n)$,
where $\psi\in\5M_k$ 
is completely unrestricted, except that it must
have a nonzero coefficient on the
monomial $z_n$. We denote the set of such
$\psi$ by $\hat{\5M_k}$, and we refer to
them as {\em admissible} elements
of ${\5M_k}$.  Thus $F$ takes the
form 
$$ F(z) = (z_1\phi_1(p,z_n),\ldots,
z_{2m}\phi_{2m}(p,z_n),\psi(p,z_n)),$$
with $\phi\in (\5F_k^\times)^{2m}$
and $\psi\in\hat{\5M_k}$.  We denote
this $F$ by $M_n(\phi,\psi)=\hat M(\phi,\psi)$.
As before, $M_n=\hat M$ is
a bijection from its domain onto $\5C_n$.

Denoting $j(z')=(z',0)$ for $z'\in\C^{2m}$,
we may write
$$ \hat M(\phi,\psi)(z',z_n) = j(z'\times\phi(p,z_{n})) + \psi(p,z_{n})e_n,$$
where $e_n$ denotes the last vector of the standard
basis of $\C^n$.

Notice that $M_n$ has a rather different
kind of domain, depending on the parity
of $n$.  

\smallskip
We shall also use the notation $J$
for the involutive element of $\5G_n$
defined by
$$ J(z)= 
\left\{
\begin{array}{rcl}
(z_2,z_1,z_4,z_3,\ldots,z_{2m},z_{2m-1})
&,& n=2m,\\
(z_2,z_1,z_4,z_3,\ldots,z_{2m},z_{2m-1},z_n)
&,& n=2m+1.
\end{array}
\right.
$$

Observe that $J$ reverses every $\L\in D_n$, i.e.\ $J^{-1}\L J = \L^{-1}$.

\subsection{The Groups $\5K_m$ and $\hat{\5K_k}$}
For $m\in\N$, we denote by $\5K_m$ the set of elements 
$F\in\5G_m$ that take the form
$$ F(t) = \left(
t_1\phi_1(t),
\ldots,
t_m\phi_m(t)
\right),
$$
with each $\phi_j\in\5F_m^\times$.
One readily checks that $\5K_m$ is a subgroup
of $\5G_m$.
\\
We use the notation $N(\phi)=N(\phi_1,\ldots,\phi_m)$
to denote $F$ of the above form.  Using the coordinatewise
product, we also write $F(t)=t\times\phi(t)$ and
$N(\phi)=\id_m\times\phi$.

\smallskip
For $k=m+1$, we denote by $\hat{\5K_k}$ the set
of elements 
$F\in\5G_k$ that (with $t=(t',t_k)$)
take the form
$$ F(t) = \left(
t_1\phi_1(t),
\ldots,
t_m\phi_m(t),
t_k\phi_k(t)+ t'\cdot\zeta(t')
\right),
$$
with each $\phi_j\in\5F_k^\times$
and $\zeta\in\5F_m$.  We remark that
{\em every} series $g(t)\in \5M_{m+1}$ with $g(0)=0$
may be written in the form
$t_k\phi_k(t)+ t'\cdot\zeta(t')$
for some $\phi_k\in \5F_k$
and $\zeta\in\5F_m$, so that the 
form of the last component $F_{,k}$ is restricted
only by the requirement that the coefficient
of the monomial $t_k$ be nonzero. This requirement
is obviously needed for the invertibility of $F$.
Thus $\hat{\5K_k}$ consists of the maps
of the form
$$ F(t) = \left(
t_1\phi_1(t),
\ldots,
t_m\phi_m(t),
\psi(t)
\right),
$$
with $\phi\in(\5F_k^\times)^m$ and
$\psi\in\hat{\5M_k}$ (i.e. $\psi$ admissible).
It is routine to check that 
$\hat{\5K_k}$ is a subgroup of $\5G_k$.

We use the notation $$\hat{N}(\phi,\psi)=\hat N(\phi_1,\ldots,\phi_k,
\psi_1,\ldots,\psi_m)$$
to denote $F$ of the above form, and we may also write
$F(t',t_k)=j(t'\times\phi(t))+\psi(t)e_k$
and
$\hat{N}(\phi,\psi)= j\circ(\id_k\times\phi) + \psi e_k$,
where $e_k$ stands for the vector $(0,\ldots,0,1)\in\C^k$.

\subsection{The Homomorphisms $P$, $H$, and $\Phi$}\label{homs}
If $n=2m=2k$ is even, then
to $F=M(\phi)\in\5C_n$ we associate the $k$ variable map
$P_n(F)\in\5G_k$ defined by
$$ P(F)(t) = t\times\pi(\phi(t))=\left(
t_1\cdot\phi_1(t)\cdot\phi_2(t),
\ldots,
t_m\cdot\phi_{2m-1}(t)\cdot\phi_{2m}(t)
\right).$$

If $n=2m+1=2k-1$ is odd, then 
to
$F=\hat M(\phi,\psi)\in\5C_n$ we associate the $k$ variable map
$$ P(F)(t) = (t',1)\times\pi(\phi(t),\psi(t))=\left(
t_1\cdot\phi_1(t)\cdot\phi_2(t),
\ldots,
t_m\cdot\phi_{2m-1}(t)\cdot\phi_{2m}(t),
\psi(t)
\right),$$
where $t=(t',t_k)$.

We have the basic semiconjugation property:
\bl For each $F\in\5C_n$, we have
\begin{equation}\label{diagram}
P(F)\circ\pi = \pi\circ F.
\end{equation}
Moreover,
property \eqref{diagram} determines $P(F)$ uniquely.
\qed
\el

Using the maps $p$ and $\pi$, we may rewrite the definition of $P$ as:
\beq
\begin{array}{rcll}
P(M(\phi)) &=& N(\pi\circ\phi)=N(p\circ\phi),&\textup{ if $n=2m$,}\\  
P(\hat M(\phi,\psi)) &=& \hat N(\pi\circ(\phi,\psi))=\hat N(p\circ\phi,\psi),&\textup{ if $n=2m+1$,}
\end{array}
\eeq 
where $\phi$ runs through $(\5F^\times_k)^{2m}$ and in the odd case
$\psi$ runs through $\hat{\5M_k}$. 

\bl $P:\5C_n\to\5G_k$ is a group homomorphism.
\el
\bpf This follows without further calculation
from the uniqueness in (\ref{diagram}) 
and the associativity of composition: If $F,G\in\5C$, then
$ P(F)P(G)\pi
=P(F)\pi G = \pi F G$, so $P(FG)=P(F)P(G)$.
\epf

However, it is useful to note the explicit formulas
for compositions of maps in the images of $M$ and $N$,
which are readily proved by direct calculation:

\bl\Label{L:M-comp}
Let $\phi$,$\phi'$ and $\phi''\in(\5F_m^\times)^{2m}$. 
Let $\chi'=P(M(\phi'))$. Then
the following are equivalent:
\\
(1)$ M(\phi'') = M(\phi)M(\phi')$.
\\
(2) $\phi''(t) = \phi'(t)\times\phi(\chi'(t))$.
\\
(3) $\phi'' =\phi'\times\phi\circ(\id_m\times\pi\circ\phi')$. 
\qed\el

\bl\Label{L:N-comp}
Let $\lambda$,$\lambda'$ and $\lambda''\in(\5F_m^\times)^m$. 
Let $\chi'=N(\lambda')$. Then
the following are equivalent:
\\
(1)$ N(\lambda'') = N(\lambda)N(\lambda')$.
\\
(2) $\lambda''(t) = \lambda'(t)\times\lambda(\chi'(t))$.
\\
(3) $\lambda'' =\lambda'\times\lambda\circ(\id_m\times\lambda')$. 
\el

Similarly, for $\hat M$ and $\hat N$, we have:

\bl\Label{L:hatM-comp}
Let $\phi$,$\phi'$ and $\phi''\in(\5F_{m+1}^\times)^{2m+1}$,
and $\psi$,$\psi'$ and $\psi''\in\hat{\5M_{m+1}}$. 
Let $\chi'=P(\hat M(\phi',\psi'))$. Then
the following are equivalent:
\\
(1)$ \hat M(\phi'',\psi'') = \hat M(\phi,\psi)\hat M(\phi',\psi')$.
\\
(2) $\phi''(t) = \phi'(t)\times\phi(\chi'(t))$ and $\psi''(t)=\psi(\chi'(t))$.
\\
(3) $\phi'' =\phi'\times\phi\circ(\id_m\times\pi\circ\phi')$
and $\psi''=\psi\circ\chi'$. 
\qed\el

\bl\Label{L:hatN-comp}
Let $\lambda$,$\lambda'$ and $\lambda''\in(\5F_m^\times)^m$
and $\psi$,$\psi'$ and $\psi''\in\hat{\5M_{m+1}}$. 
Let $\chi'=\hat N(\phi',\psi')$. Then
the following are equivalent:
\\
(1)$ \hat N(\lambda'',\psi'') = \hat N(\lambda,\psi)\hat N(\lambda',\psi')$.
\\
(2) $\lambda''(t) = \lambda'(t)\times\lambda(\chi'(t))$
and $\psi''(t)=\psi(\chi'(t))$.
\\
(3) $\lambda'' =\lambda'\times\lambda\circ(\id_m\times\lambda')$
and $\psi''=\psi\circ\chi'$. 
\el

The fact that $P$ is a homomorphism is obtained again in the even case
by precomposing $\pi$ with
the equation in part (3) of Lemma \ref{L:M-comp}
and using $P(M(\phi))=\id_m\times\pi\circ\phi$. In fact,
$\pi(a\times b)=\pi(a)\times\pi(b)$, so (3) gives
 $$\pi\circ\phi'' =(\pi\circ\phi')\times(\pi\circ\phi)\circ(\id_m\times\pi\circ\phi').$$ 

A similar argument applies in the odd case, using Lemmas
\ref{L:hatM-comp} and \ref{L:hatN-comp}.

\medskip
The kernel of $P$
is the set of maps $F$ of the form
$$ F(z) = 
\left\{
\begin{array}{rcl}
M(\phi)&,& n=2m=2k,\\
M(\phi,1)&,& n=2m+1=2k-1
\end{array}
\right.
$$
where $\phi\in(\5F_k^\times)^{2m}$
and $\phi_{2j-1}(t)=1/\phi_{2j}^{-1}(t)$
for $j=1,\ldots,m$.
The group $\ker P$ is abelian.

\medskip

For $n=2m$, the map
$$\Phi:
\left\{
\begin{array}{rcl}
(\5F_k^\times)^m &\to& \5C_n,\\
\phi &\mapsto& 
M\left(\phi_1, 1/\phi_1,\ldots,\phi_m,
1/\phi_m)
\right),
\end{array}
\right.
$$
is a group isomomorphism onto $\ker P$.

For $n=2m+1$, the corresponding isomorphism onto $\ker P$
is 
$$\Phi:
\left\{
\begin{array}{rcl}
(\5F_k^\times)^m &\to& \5C_n,\\
\phi &\mapsto& 
M\left(\phi_1, 1/\phi_1,\ldots,\phi_m,
1/\phi_m,1)
\right),
\end{array}
\right.
$$

\medskip
Note that in each case the image of $\Phi$
consists of reversible elements. All are reversed by
$J$.
\sr{added note}

\medskip
Clearly, when $n=2m$, 
the image of $P$ lies in $\5K_m$,
and when $n=2m+1$, the image lies in
$\hat{\5K_k}$.
To see that these are the exact images
of $P$ in the respective cases, we
define right inverse maps:

For $n=2m$, we define
$H:\5K_m\to\5C_n$ by 
\beq\Label{H-even}
H(N(\lambda))(z) = M(\epsilon\circ\lambda)=
M(\lambda_1,1,\lambda_2,1,\ldots,\lambda_m,1).
\eeq
From the definition of $P$, we see that
$P(M(\phi)) = N(\pi\circ\phi)$
so obviously $P(H(N(\lambda)))=N(\pi\epsilon\lambda)=N(\lambda)$, as required.

(Moreover,
for each $\chi\in\5K_m$,
$H(\chi)$ is the unique element $F\in\5C_n$ with the properties
$F(z)_{,j}=z_j$ for all even $j$ and 
$F\circ\epsilon=\epsilon\circ\chi$.)

For $n=2k-1$, we define
$H=H_n:\hat{\5K_k}\to \5C_n$ by
\beq\Label{H-odd}
H(\hat N(\lambda,\psi))(z) = 
\hat M(\epsilon\circ(\lambda,\psi))(z)=\hat M(\lambda_1,1,\lambda_2,1,\ldots,\lambda_m,1,\psi)(z).
\eeq
In this case, the definition of $P$ amounts to
$$
P(\hat M(\phi,\psi)) = 
\hat N(p\circ\phi,\psi),
$$
so again $P(H(N(\lambda,\psi)))=N(\lambda,\psi)$, as required.

(Moreover,
for each $\chi\in\hat{\5K_k}$,
$H_n(\chi)$ is the unique element $F\in\5C_n$ with the properties
$F(z)_{,j}=z_j$ for all even $j$ and 
$F\circ\epsilon=\epsilon\circ\chi$.)

So we have proved:

\bl The image of $P$ is $\5K_m$ if $n=2m$
and is $\hat{\5K_k}$ if $n=2k-1$. 
\qed\el

Next we have:
\bl For each $n\ge2$, $H_n$ is a group homomorphism.
\el
\bpf 
We give the explicit version, taking the cases separately.

\medskip
\noindent
{$1^\circ$}: $n=2m$. 
Fix two elements $\chi$,$\chi'\in\5K_m$, and let $\chi''=\chi\chi'$. 
There are unique elements $\lambda$,$\lambda'$ and 
$\lambda''\in(\5F_m^\times)^m$ 
such that
$\chi=N(\lambda)$,
$\chi'=N(\lambda')$,
and $\chi''=N(\lambda'')$.
By Lemma \ref{L:M-comp} and the definition of $H$, 
it suffices to show that
$$\epsilon\circ\lambda'' =\epsilon\circ\lambda'\times
\epsilon\circ\phi\circ(\id_m\times\pi\circ\epsilon\circ\phi')$$
But this is immediate from Lemma \ref{L:N-comp}
and the fact that $\pi\circ\epsilon=\id_m$.

\medskip
\noindent
{$2^\circ$}: $n=2m+1$.
Fix $\chi$,$\chi'\in\hat{\5K_k}$,
and let $\chi''=\chi\chi'$.
There are unique elements $\lambda$,$\lambda'$ and 
$\lambda''\in(\5F_k^\times)^m$ 
and $\psi$,$\psi'$ and $\psi''\in\hat{\5M_k}$
such that
$\chi=\hat N(\lambda, \psi)$,
$\chi'=\hat N(\lambda',\psi')$,
and $\chi''=\hat N(\lambda'',\psi'')$.
By the definition of $H$, we have to show that
$$ \hat M(\epsilon\circ\lambda'',\psi'')=
\hat M(\epsilon\circ\lambda,\psi)\hat M(\epsilon\circ\lambda',\psi').$$
By Lemma \ref{L:hatM-comp}, in view of the fact that
$$P(\hat M(\epsilon\circ\lambda',\psi'))=
\hat N(\pi\circ\epsilon\circ(\lambda',\psi'))
=\hat N(\lambda',\psi')=\chi',$$
this amounts to
showing that
$$ \epsilon\circ\lambda''=
\epsilon\circ\lambda'\times\epsilon\circ\lambda\circ\chi' 
\textup{ and }
\psi''=\psi\circ\chi'.$$
But this is immediate from Lemma \ref{L:hatN-comp}.
\epf

\begin{figure*}
\vskip3mm
$$
\begin{array}{rcccccccl}
 (1) &\to& (\5 F_m^\times)^m
&\xrightarrow[\Phi]{}&  
\5 C_{2m}
& 
\overset{\xleftarrow{H}}{\xrightarrow[P]{}} 
&\5 K_m
&\to(1)\\
&&&&\cap&&\cap&&\\
&&&&\mymf G_{2m}&&\5 G_m&&\\
&&&&&&&&\\
&&&&&&&&\\
(1) &\to& (\5 F_{m+1}^\times)^m
&\xrightarrow[\Phi]{}&  
\5 C_{2m+1} 
&
\overset{\xleftarrow{H}}{\xrightarrow[P]{}} 
&
\hat{\5 K}_{m+1}
&\to&
(1)\\
&&&&\cap&&\cap&&\\
&&&&\5 G_{2m+1}&&\5 G_{m+1}&&\\
\end{array}
$$
\vskip3mm
\caption{Exact Sequences}
\label{fig:exact-sequence2}
\end{figure*}

\bc\Label{C:semidirect}
Let $n\ge2$. The group
$\5C_n$ is the semidirect
product of $\im\Phi$ and $\im H$.
Each $F\in\5C_n$ has a unique factorization
in the form 
$H(\chi)\Phi(\phi)$, 
with $\chi\in\im P$ and
$\phi\in(\5F_k^\times)^m$. 
\ec
\bpf This follows from the facts that
$\im\Phi=\ker P$ is normal and that
the homomorphism $H$ is a right inverse for $P$.
\epf
However, $\5C_n$ is not
the {\em direct} product of $\im H$ and $\im\Phi$. 

The structural results of this subsection are
\sr{added remark about exact sequences}
summarized in Figure \ref{fig:exact-sequence2}

\section{Proof of Theorem \ref{T:main} in Dimension $2$}\Label{S:D-2}
\sr{modified text}
The case $n=2$ of our main theorem serves as the
foundation layer for an inductive proof of the general
case, and now we lay this down.

For the reader's convenience,
We include the short 
proof of our previously-published 
Theorem \ref{T:prod-4-rev} (which is the same as
part (1) of Theorem \ref{T:main} in case $n=2$). 

\bpf
Let $F\in\5G_2$ have $\det L(F)=1$.
We have to show that $F$
it may be factorized as $F=
g_1 g_2 g_3 g_4$, where each $g_j$
is reversible in $\5G$.

In fact, if $\det L(F) =1$, then
multiplying by some (reversible) $\Lambda\in D$ (possibly the identity)
we can arrange that $L(F\Lambda)$ is conjugate to an
infinite-order element of $D$. Then by Poincar\'e-Dulac, $F\Lambda$
is conjugate (say by $K\in\5G$) 
to some element of $\5C$, so $(F\Lambda)^K$ may be
factored as $H(\chi) \Phi(\phi)$, where
$\chi(t)=t+ \HOT$. Now $\Phi(\phi)$
is reversible, 
and we know \cite[Theorem 9(2)]{O} 
that $\chi$ is the product
of two reversibles in $\5G_1$, so $H(\chi)$ is the
product of two reversibles, say $H(\chi_1)$ and $H(\chi_2)$.
Thus 
$$F^K= H(\chi_1) H(\chi_2) \Phi(\phi)
(\Lambda^{-1})^K$$
is the product of four reversibles, and conjugating with $K^{-1}$ we obtain the result.
\epf

\medskip
\bpf[Proof of Theorem \ref{T:main} part (2) when $n=2$]
Let  $F\in\5G_2$  with $\det L(F)=+1$. 
With the notation
in the last proof,
$\Lambda^K$ 
and $\Phi(\phi)$ are strongly-reversible,
since $\Lambda$ and $\Phi(\phi)$ are
reversed by the involution $J$.
Thus it suffices to prove that $H(\chi)$ 
is the product of $10$ involutions, whenever
$\chi(t)=t+\HOT\in\5G_1$.

Now we may factor $\chi$ as $\chi_1\chi_2$, where
these take the form
$$\begin{array}{rcl}
\chi_1(t) &=& t+\alpha t^2+\alpha^2 t^3,\\
\chi_2(t) &=& t + \beta t^3 +\HOT.
\end{array}
$$
(possibly with $\alpha=0$ or $\beta=0$).  
The map
$\chi_1$ is the identity or
is conjugate to the map
$f_1$, given by (cf. \cite{O})
\beq\Label{1D-normal} f_1(t) = \frac{t}{1-t}
=
t+t^2+t^3+t^4+\HOT,
\eeq
reversed by $t\mapsto -t$,
hence is
the product of two involutions in $\5G_1$,
and hence so is $H(\chi_1)$.  Thus, since
$H(\chi)=H(\chi_1)H(\chi_2)$, it suffices to show
that $S=H(\chi_2)$ is the product of $8$ involutions.

If $\beta=0$, then $\chi_2$ is the product of 
$4$ involutions in view of \cite[Theorem~9]{O}, so we consider the case $\beta\not=0$.
By a conjugation in $\5G_1$, we may take $\beta=1$, so
$S$ takes the form
$$ S(z) = \left( z_1( 1+p^2+\HOT), z_2\right).$$
Take the map
$$
T=\left(\begin{matrix} a&b\\c&d \end{matrix}\right)\in\SL(2,\C) .  
$$
A calculation yields that, up to terms of degree 5,
$$T^{-1}ST(z) =
\begin{pmatrix}
d&-b\\
-c&a
\end{pmatrix}
\begin{pmatrix}
1+ A&0\\
0&1
\end{pmatrix}
\begin{pmatrix}
az_{1}+bz_{2}\\
cz_{1}+dz_{2}
\end{pmatrix}
=
\begin{pmatrix}
z_1+dB\\
z_2- cB
\end{pmatrix},
$$
where $\det T=1$, $A=(az_1+bz_2)^2(cz_1+dz_2)^2$ and $B=(az_1+bz_2)^3(cz_1+dz_2)^2$. 

Let $\Lambda_1 = \left(\begin{matrix}2&0\\0&\frac12\end{matrix}\right)$. 
Then $\L_1 T^{-1}ST$ equals
\begin{equation}\Label{res}
\left(2z_{1}\left(1+{ad}
(a^{2}d^{2}+6abcd+3b^{2}c^{2}) p^{2}\right), 
\frac12 z_{2} \left(
1-{bc}(b^{2}c^{2}+6 abcd + 3a^{2}d^{2})
p^{2}\right)
\right)
\end{equation}
plus non-resonant terms of order $5$.
By the Poincar\'e-Dulac Theorem, $\L_1 T^{-1}ST$ is
conjugate to the map $S_1$ obtained by removing all
non-resonant terms.  A calculation shows
that $S_{1}$ equals \eqref{res} up to terms of degree $5$ in $z$.
We now choose $T$ such that
\begin{equation}
ad(a^{2}d^{2}+6abcd+3b^{2}c^{2})=bc(b^{2}c^{2}+6 abcd + 3a^{2}d^{2}),
\end{equation}
or, substituting $ad=bc+1$,
\begin{equation}
(bc+1)((bc+1)^{2}+6bc(bc+1)+3b^{2}c^{2})=bc(b^{2}c^{2}+
6bc(bc+1) + 3(bc+1)^{2}),
\end{equation}
which simplifies to
\begin{equation}
6b^{2}c^{2} + 6 bc +1
 =0
\end{equation}
and clearly has a solution.

Then $S_1(z)$  factors as
$H(\chi)\Phi(\phi)$ with $\chi(t)=P(S_{1})(t)=t + O(t^{4})$.
Hence by \cite[Theorem~9]{O}, $\chi$ and therefore $H(\chi)$ is the product
of four involutions. Thus $S_1$
is the product of $6$ involutions.  
Thus $\L_1 T^{-1}ST$ is the product of $6$
involutions, so $S$ is the product of $8$.
This concludes the proof.
\epf

Each product 
$F=f_1\cdots f_n$
of reversible $f_j$'s has 
$\det L(F)=\pm1$, so (multiplying if
necessary by a suitable linear involution)
it follows from Theorem \ref{T:prod-4-rev} that each product of
reversibles reduces to the product of
five.  It also follows that the elements that
are products of reversibles are precisely those
with $\det L(F)=\pm1$.
\hl{Thus the case $n=2$ of Corollary \ref{C:main} is immediate.}

\section{Proof of Theorem \ref{T:main} in Dimension $n>2$}\Label{S:D-3}

We will actually prove a more refined result,
in which the number of factors required depends
in a more complicated way on the dimension $n$.

First, we introduce notation for the number of factors needed,
in various situations:

For $n\ge2$,  
let $r_1(n)$ denote the least
$r\in\N$ such that each $F\in\5G_n$ having
$\det L(F)=1$ may be expressed as the product
of $r$ reversible elements of $\5G_n$.
Similarly, let $r_d(n)$ be the least number
of reversible factors from $\5G_n$ required for the factorization
of each $F\in\5G_n$ having $L(F)\in D_n$.
Finally, let $r_c(n)$ be the least number
of reversible factors from $\5G_n$ required for the factorization
of each $F\in\5C_n$ having $L(F)\in D_n$.

It is obvious that
\beq\Label{E:obvious-r}
r_c(n) \le r_d(n)\le r_1(n)
\eeq
whenever $n\ge2$.

\bl\Label{L:factor-diag-2} Let $n\ge2$.
Each diagonal matrix $T\in\SL(n,\C)$
may be factored as the product
of two diagonal matrices $T_1T_2$,
where $T_1\in D_n$ and 
there is a permutation
matrix $\sigma$ such that
the conjugate $T_2^\sigma$
belongs to $D_n$.
\el

\bpf
Let $T=\diag(\l_1,\l_2,\ldots,\l_n)$,
and note that
$\l_1\cdots\l_n=1$.

Take 
$$T_1=\diag(\l_1,1/\l_1,\l_1\l_2\l_3, 1/(\l_1\l_2\l_3),\ldots),$$
and $$T_2=\diag(1,\l_1\l_2,1/(\l_1\l_2),\l_1\l_2\l_3\l_4,
1/(\l_1\l_2\l_3\l_4),\ldots).$$

If $n$ is odd, 
then the last entry in $T_1$ is $1$,
so $T_1\in D_n$. Since $T_2$ is conjugated into $D_n$
by the permutation $(1n)$ that swaps the coordinates $z_1$ and $z_n$, 
we are done, in this case.

If $n$ is even, then $T_1\in D_n$, and $T_2$ has both
first and last entries equal to $1$, so it is 
conjugated into $D_n$ by the $n$-cycle $(12\ldots n)$
that rotates the last coordinate back into first
position, and shifts the others down.
\epf

\bl\Label{L:factor-Dn-2}
Each element of $D_n$ may be factored as the product
of two generic elements of $D_n$.
\el

\bpf  Let $T\in D_n$. Then $T=\Phi(\alpha)$
for some $\alpha\in(\C^\times)^m$. Choose $\l\in(\C^\times)^m$
such that $\l_j$ is \hl{multiplicatively independent of
$\alpha_1$,$\ldots$,$\alpha_m$,$\l_1$,$\ldots$,$\l_{j-1}$, for each $j$.}
Take $T_1=\Phi(\alpha\times\lambda)$
and $T_2=\Phi(\lambda_1^{-1},\ldots,\l_m^{-1})$. Then
each $T_j$ is a generic element of $D_n$, and $T=T_1T_2$.
\epf

\bl\Label{L:factor-diag-3} Let $n\ge2$.
Each diagonal matrix $T\in\SL(n,\C)$
may be factored as the product
of three diagonal matrices  $T_1T_2T_3$,
where $T_1\in D_n$ and 
there is a permutation
matrix $\sigma$ such that
the conjugates $T_2^\sigma$
and $T_3^\sigma$
belong to $D_n$, and $T_3^\sigma$
is a generic element. 
\el

\bpf 
Let $T_1T_2'$ be the factorization and $\sigma$ the permutation
given by Lemma \ref{L:factor-diag-2},
and apply Lemma \ref{L:factor-Dn-2} to $(T_2')^\sigma$.
\epf

\bl\Label{L:to-generic} 
Let $n\ge3$. Then $r_1(n)\le r_c(n)+2$.
\el

\bpf
Fix $F\in\5G_n$
with $\det L(F) =1$. 

By using a linear conjugation, if
need be, we 
may assume that $L(F)$ is in Jordan
canonical form, so that the
diagonal elements multiply to $1$.

Write $L(F)=T+N$,
where $T$ is diagonal and $N$
is strictly upper triangular.
Applying the last lemma,
we can write $T=T_1T_2T_3$,
where $T_1\in D_n$ and 
both $T_2$ and $T_3$ are diagonal, and conjugate
by the same permutation $\sigma$ of coordinates
to elements of $D_n$, with $T_3^\sigma$ 
generic.
Let $F_1=(T_1T_2)^{-1}F$.
Then $L(F_1)$ is upper triangular,
with the same diagonal  as $T_3$.

The eigenvalues of $L(F_1)$ are its diagonal elements,
and are distinct, so we may conjugate $L(F_1)$ to $T_3$
by using an element of $\GL(n,\C)$. Applying the
same conjugation to $F_1$, we conjugate $F_1$ to a map
$F_2$ with $L(F_2)=T_3$.  Applying Poincar\'e-Dulac, we can conjugate $F_2$ 
to a map $F_3$ that commutes with $T_3$, without changing
the linear part, so $L(F_3)=T_3$.  Then
$F_3^\sigma$ commutes with $T_3^\sigma$,
and hence belongs to $\5C_n$, and has $L(F_3^\sigma)=T_3^\sigma$.  

Now $F_3^\sigma$
is the product of $r_c(n)$ 
reversibles, hence so are $F_3$, $F_2$
and $F_1$. Since $T_1$ and $T_2$ are reversible,
$F$ is the product of $2+r_c(n)$ reversibles. 
\epf	

\bl\Label{L:D-to-generic} 
Let $n\ge3$. Then $r_d(n)\le r_c(n)+1$.
\el

\bpf
Fix $F\in\5G_n$ with $L(F)\in D_n$. By Lemma \ref{L:factor-Dn-2},
we may factor $L(F)=T_1T_2$, where each $T_j\in D_n$ is generic.
Taking $F_1=T_1^{-1}F$, we have $L(F_1)=T_2$, and applying
Poincar\'e-Dulac we can conjugate $F_1$ to an element
$F_2$ of $\5C_n$ having $L(F_2)=T_2$. Since
$F_2$ is the product of $r_c(n)$ reversibles, so is $F_1$,
and hence $F=T_1F_1$ is the product of $1+r_c(n)$.
\epf

\bl\Label{L:to-half}
Let $m\ge1$. Then\\
(1) $ r_c(2m)\le 1+ r_d(m)$, and\\
(2) $r_c(2m+1)\le 1+r_1(m+1)$.
\el

\bpf (1) Let $n=2m$.
Fix $F\in\5C_n$, with $L(F)\in D_n$.
Then
$\chi=P(F_1)$ belongs to $\5G_m$ and
is tangent to the identity, 
so it may be factored as the product
of $r_d(m)$ reversibles.

By Corollary \ref{C:semidirect}, 
we can factor $F$ as $H(\chi)\Phi(\phi)$,
for some $\phi\in(\5F_m^\times)^m$,
and we know that $\Phi(\phi)$ is reversed by $J$,
so $F$ is the product
of $1+r_d(m)$ reversibles.
Thus $r_c(n)\le 1+r_d(m)$.

\smallskip
\noindent(2)
Let $n=2m+1=2k-1$.
Fix $F\in\5C_n$, with $L(F)\in D_n$.
Then this time
$\chi=P(F_1)\in\5G_k$ may fail to be tangent to the identity,
or even to belong to $D_k$,
but still has $\det L(\chi)=1$,
so it may be factored as the product
of $r_1(k)$ reversibles.
Proceeding as before, we get
$r_c(n)\le 1+r_1(k)$,
as required.
\epf	

\bc\Label{C:r1-case} If $n\ge2$, 
then $r_1(n)\le 1+3\cdot\thec$.
\ec

\bpf  
We proceed inductively, starting at $n=2$.

For $n=2$, Theorem \ref{T:prod-4-rev} tells us that
$r_1(n)\le4=1+3\cdot\thec$.

Fix $n>2$, and assume that 
for every $n'<n$,
we have $r_1(n')\le 1+3\cdot\textup{ceiling}(\log_2 n')$.

Then with $k$ as usual, Lemmas \ref{L:to-generic}
and \ref{L:to-half} and
inequalities \ref{E:obvious-r} yield
$$ r_1(n) \le 2+r_c(n)\le 3+r_1(k)
\le 4+3\cdot\textup{ceiling}(\log_2 k) $$
so, since $\textup{ceiling}(\log_2 k)$ 
is one less than $\thec$, we have
$r_1(n)\le1+3\cdot\thec$, and the induction step
is complete.
\epf

This Corollary has the same content as Theorem \ref{T:main}, part (1),
so that is now proven.

\bpf[Proof of Theorem \ref{T:main}, Part (2)]
Denote the minimal number of involutive factors needed
to express each 
member of the classes corresponding to
$r_1$, $r_d$ and $r_c$, respectively, by $i_1$,
$i_d$ and $i_c$, respectively. Observing that
the elements of $D_n$ and of $\im\Phi$ are strongly-reversible, 
and reviewing the proofs of 
Lemmas \ref{L:to-generic}
and \ref{L:to-half},
we obtain the following estimates:
$$\begin{array}{rcl}
i_c(n)&\le& i_d(n)\le i_1(n)\\
i_1(n)&\le& i_c(n)+4\\
i_d(n)&\le&i_c(n)+2\\
i_c(2m)&\le&2+i_d(m)\\
i_c(2m+1)&\le&2+i_1(m+1),
\end{array}
$$
whenever $m,n\in\N$ and the terms on both sides are defined (i.e.
we say nothing about $i_1(1)$, $i_d(1)$
or $i_c(1)\,)$.
We can now carry out an induction
to estimate $c_1(n)$, and each induction step
adds $6$ to the number of involutions that will
suffice.

At the lowest level, when $n=2$, Theorem \ref{T:prod-4-rev} part (2)
tells us that $14=\answerIS$ involutions suffice,
so induction gives the
result, since $\thec$ increases by $1$ at each step.
\epf

\bpf[Proof of Corollary \ref{C:main}]
The equivalence of (1), (2) and (3)
follows from the theorem and the fact
that each reversible, and hence each product of reversibles 
has determinant $\pm1$.

Closer analysis of the proof of the theorem given above reveals
that each $F$ with $\det L(F)=1$ may also be represented
\sr{referee says it's heavily convoluted. I'm inclined to leavit alone.}
as the product of $6\cdot\thec$ involutions
and one special map
that is a homomorphic
image of an element $\chi\in\5G_1$ having
multiplier $+1$. (The homomorphism is the composition
of repeated $H$ and inner automorphisms.)
Examining the detail
in the proof of Theorem \ref{T:O},
one finds that $\chi$ is the product
of two reversibles, one strongly reversible,
and the other reversed by an element
of order dividing $4$.
(The theorem is Theorem 9 of [O],
and the proof is on pp. 18-19 of that paper.
The map is denoted $f$, instead of $\chi$.
Three cases are considered. In case $1^\circ$,
$f$ is factored as $gh$, where
$g$ is conjugate to $z + z^2 + z^3$, which is
strongly reversible, and $h$ is $\id$
or is conjugate to $z+z^3+\frac32 z^5$, which is reversed
by $z\mapsto iz$. In case $2^\circ$ --- note that
there is a misprint: this case is $p>2$, not
$p\ge2$ ---, $f$ is the product of two maps
conjugate to $z+z^2+z^3$. Finally, in case $3^\circ$,
$f=gh$, where $g$ is conjugate to 
the aforementioned $z+z^3+\frac32 z^5$
and $h$ is conjugate to
$z+z^4+2z^7$, and hence is strongly reversible.) 
Thus $\chi$, and hence the special map,
are each the product of two involutions
and two reversible maps of degree dividing $4$,
\sr{footnote removed, and used in 
introduction}
so that
$F$ is the product of $2+6\cdot\thec$ involutions
and two reversible maps of degree dividing $4$.
\epf

\begin{table}[ht!]
\begin{center}
\begin{tabular}{|c|c|c|c||c|c|c|c|}
\hline
$n$& $r_1(n)$ & $r_d(n)$&$r_c(n)$&$n$& $r_1(n)$ & $r_d(n)$&$r_c(n)$\\
\hline
\hline
1&2&2&2 &9&13&12&11\\
2&4&4&3 &10&12&11&10\\
3&7&6&5 &11&12&11&10\\
4&7&6&5 &12&11&10&9\\
5&10&9&8&13&13&12&11\\
 
6&9&8&7 &14&12&11&10\\
7&10&9&8&15&12&11&10\\
 
8&9&8&7& 16&11&10&9\\
\hline
\end{tabular}
\end{center}
\caption{}\Label{Table}
\end{table}

\br The inequalities in Lemmas
\ref{L:to-generic} and \ref{L:to-half}
may be used to derive
estimates for $r_1(n)$ that are often considerably
smaller than the estimate $\answerRS$.
These estimates depend on the parity
of the terms in the chain of links $n'\to k'$
connecting $n$ to $2$.
For instance, from the chain
$$96\to48\to24\to12\to6\to3\to2$$
one obtains $r_1(96)\le14$,
in contrast to the estimate $r_1(97)\le 20$
obtained from the chain
$$ 97\to49\to25\to13\to7\to4\to2.$$  
The best estimates
are obtained for powers of $2$: 
$$ c_1(2^n) \le 2+2n.$$
Table \ref{Table} gives the best estimates
obtainable from these Lemmas for the first few 
$n$.

We do not know sharp values for $r_1(n)$ or $r_d(n)$, in any case
of dimension greater than $1$.  
\er


\begin{thebibliography}{99}

\bibitem{AO} {\bf P. Ahern and A.G. O'Farrell.}
Reversible biholomorphic germs.  
{\em Comput. Methods Funct. Theory} {\bf 9} (2009), 473--84.

\bibitem{AG} {\bf P. Ahern and X. Gong.} A Complete Classification for Pairs of Real Analytic Curves in
the Complex Plane with Tangential Intersection. {\em J. Dynamical and Control Systems} {\bf 11}
(2005), 1-71.

\bibitem{B} {\bf G.D. Birkhoff.} The restricted problem of three bodies, 
{\em Rend. Circ. Mat. Palermo} {\bf 39} (1915), 265--334.

\bibitem{BZ} {\bf F. Bracci and D. Zaitsev.} Dynamics of one-resonant biholomorphisms.
{\em Journal of the European Mathematical Society}, to appear.
\sr{appeared?}

\bibitem{BrFa}
{\bf T.E. Brendle and B. Farb.}
{Every mapping class group is generated by 6 involutions.}
{\em J. Algebra}
{\bf 278}
(2004)
{187--98}.

\bibitem{C} {\bf R. Camina.} The Nottingham Group. 
{\em New Horizons in pro-$p$ groups. Progr. Math.} 
{\bf 184} (2000) 205-21.
  
\bibitem{Dj1967}
{\bf D.{\v{Z}}. Djokovi{\'c}.}
{Product of two involutions.}
{\em Arch. Math. (Basel)}
{\bf 18}
(1967)
{582--84}.

\bibitem {Dj1986}
{\bf D.{\v{Z}}. Djokovi{\'c}.}
{Pairs of involutions in the general linear group.}
{\em J. Algebra.}
{\bf 100}
(1986)
{214--23.}

\bibitem{El1993}
{\bf E.W. Ellers.}
{The reflection length of a transformation in the unitary group
              over a finite field.}
{\em Linear and Multilinear Algebra}
{\bf 35}
(1993)
{11--35}.

\bibitem{FZ}
{\bf W. Feit and G.J. Zuckerman.}
{Reality properties of conjugacy classes 
in spin groups and symplectic groups.}
{\em Contemp. Math.}
{\bf 13}
(1982)
{239--53}.

\bibitem{FS} {\bf N. J. Fine and G. E. Schweigert.}
On the group of homeomorphisms of an arc. {\em Ann. of Math.}
{\bf 62} (2) (1955) 237--53.

\bibitem{GS}
{\bf {\'E}. Ghys  and V. Sergiescu.}
{Stabilit\'e et conjugaison diff\'erentiable pour certains
              feuilletages}.
{\em Topology}
{\bf 19}
(1980)
{179--97}.

\bibitem{G} {\bf X. Gong.} Anti-holomorphically reversible holomorphic 
maps that are not holomorphically reversible. 
Geometric function theory in several complex variables, 
151--64, World Sci. Publishing, River Edge, NJ, 2004.

\bibitem{Go76}
{\bf R. Gow.}
{Real-valued characters and the {S}chur index}.
{\em J. Algebra.}
{\bf 40}
(1976)
{258--70.}

\bibitem{Go79}
{\bf R. Gow.}
{Real-valued and {$2$}-rational group characters.}
{\em J. Algebra}.
{\bf 61}
(1979)
{2}
{388--413}.

\bibitem{Go81}
{\bf R. Gow.}
{Products of two involutions in classical groups of
              characteristic {$2$}.}
{\em J. Algebra}
{\bf 71}
(1981)
{583--91}.

\bibitem{GHR}
{\bf W.H. Gustafson and P.R. Halmos,  and H. Radjavi.}
{Products of involutions}.
{\em Linear Algebra and Appl.}
{\bf 13}
(1976)
{157--62}.

\bibitem{PD} {\bf Yu. Ilyashenko and S. Yakovenko. }
Lectures on Analytic Differential Equations.
AMS. Providence. 2008. (ISBN: 978-0-8218-3667-5).

\bibitem{Is}
{\bf H. Ishibashi.}
{Involutary expressions for elements in {${\rm GL}\sb n({\bf
              Z})$} and {${\rm SL}\sb n({\bf Z})$}.}
{\em Linear Algebra Appl.}
{\bf 219}
(1995)
{165--77}.

\bibitem{K} {\bf B. Klopsch.} Automorphisms of the Nottingham group.
{\em J. Algebra} {\bf 223} (2000), 37-56.
  
\bibitem{KL} {\bf B. Klopsch and V. Lev.} How long does it take
to generate a group? {\em J. Algebra} {\bf 261} (2003), 145-71.

\bibitem{KnNi}
{\bf F. Kn{\"u}ppel and K. Nielsen.}
{Products of involutions in {${\rm O}\sp +(V)$}}.
{\em Linear Algebra Appl.}
{\bf 94}
(1987)
{217--22}.

\bibitem{MW} {\bf J.K. Moser and S.M. Webster.} Normal forms for real 
surfaces in ${\C}\sp{2}$ near complex tangents and hyperbolic 
surface transformations. {
\em Acta Math.} {\bf 150} (1983), no. 3-4, 255--96.

\bibitem{Ni1987}
{\bf K. Nielsen.}
{On bireflectionality and trireflectionality of orthogonal
              groups.}
{\em Linear Algebra Appl.}
{\bf 94}
 (1987)
{197--208}.
 
\bibitem{O} {\bf A.G. O'Farrell.}
Compositions of involutive power series and reversible series. 
{\em Comput. Methods Funct. Theory} {\bf 8} (2008) 173-93. 

\bibitem{OCRM} {\bf A.G. O'Farrell.}
Reversibility questions in groups arising in analysis. 
{\em CRM Proceedings and Lecture Notes} {\bf 55} (2012) 293-300.

\bibitem{OZ} {\bf A.G. O'Farrell and D. Zaitsev.} Formally reversible 
maps of $(\C^2,0)$. {\em Ann. SNS Pisa}, to appear. 
\begin{verbatim}DOI: 10.2422/2036-2145.201201_001\end{verbatim}

\bibitem{Sar}
{\bf P. Sarnak.}
Reciprocal geodesics.
{\em Clay Math. Proc.}
{\bf 7}
(2007)
217--37.
\bibitem{St}
{\bf A. Stein.}
{{$1\frac12$}-generation of finite simple groups.}
{\em Beitr\"age Algebra Geom.}
{\bf 39} 
(1998)
{349--58}.

\bibitem{W} {\bf S.M. Webster} Pairs of intersecting real manifolds in complex space. {\em Asian J. Math.} 
{\bf 7} (2003), no. 4, 449-62. 

\end{thebibliography}
\end{document}